\documentclass[11pt,reqno]{amsart}

\usepackage{geometry}
\usepackage{amsmath}
\usepackage{amsfonts}
\usepackage{amssymb}
\usepackage{amsthm}
\usepackage{amscd}
\usepackage{stmaryrd,mathrsfs}
\usepackage{enumerate}
\usepackage{upgreek}
\usepackage{fancyhdr}
\usepackage{hyperref}

\newtheorem{theorem}{Theorem}[section]
\newtheorem{corollary}[theorem]{Corollary}

\newtheorem{remark}[theorem]{Remark}
\DeclareMathOperator{\card}{card}

\makeatletter
\def\imod#1{\allowbreak\mkern6mu({\operator@font mod}\,\,#1)}
\makeatother

\begin{document}
\large
\title{Comment on ``A note on the squarefree density of  polynomials}
\author{Yuri G. Zarhin}

\address{Department of Mathematics, Pennsylvania State University, University Park, PA 16802, USA.}
\email{zarhin@math.psu.edu}


\begin{abstract}
This is an exposition of results of a joint paper of R.C. Vaughan and the author (Mathematika
{\bf 70} (2024), no. 4).
We discuss how often the  squarefree values  of an integral polynomial  do occur.
\end{abstract}

\maketitle

As usual, $\mathbb{Z}$ and $\mathbb{Q}$ stand for the ring of integers and the field of rational numbers respectively. If $s$ is a positive integer then  $\mathbb{Z}[X_1, \dots, X_s]$ and $ \mathbb{Q}[X_1, \dots, X_s]$ stand for the rings of polynomials in $s$ variables over $\mathbb{Z}$ and $\mathbb{Q}$ respectively.
If $A$ is a finite set then we write $\card(A)$ for its cardinality.

Lel $\mathcal{P}(X_1, \dots, X_s)\in \mathbb{Z}[X_1, \dots, X_s]$ be a polynomial of degree $d \ge 2$  with integer coefficients. If $p$ is a prime then let us put
$$ \rho_{\mathcal P}(p^2): = \card\{(a_1, \dots, a_s)\in \mathbb{Z}^s/p^2 \mathbb{Z}^s \mid \mathcal P(a_1, \dots, a_s)\equiv 0\imod p^2\}.$$
 Clearly,
 $$0 \le \rho_{\mathcal P}(p^2)\le p^{2s}, \quad 0 \le  1-\frac{\rho_{\mathcal P}(p^2)}{p^{2s}} \le 1.$$

 {\bf Definition}.
 We say that
 $\mathcal{P}$ enjoys property (a) if there is a prime $p$  such that 
 $$\mathcal{P}(a_1,\dots, a_s) \in p^2\mathbb{Z} \ \ \forall  (a_1,\dots, a_s) \in \mathbb{Z}^s,$$
 i.e.,  there exists a prime $p$ such that
 $$ \rho_{\mathcal P}(p^2)=p^{2s}, \quad 1-\frac{\rho_{\mathcal P}(p^2)}{p^{2s}}=0.$$
\vskip .2cm

 Let us consider the corresponding Euler product
 $$\mathfrak S_{\mathcal P} = \prod_p \left(
1-\frac{\rho_{\mathcal P}(p^2)}{p^{2s}}
\right)$$
where $p$ runs through the set of all primes. Clearly, $\mathfrak S_{\mathcal P}=0$ if 
${\mathcal P}$ enjoys property (a).
\vskip .2cm

{\bf Definition}.
We say that $\mathcal{P}$ is {\sl squarefree} if it enjoys the following {\sl equivalent} properties.
\begin{itemize}
\item[(i)]
There is a non-constant polynomial $\mathcal{L}$ in $\mathbb{Z}[X_1, \dots, X_s]$ such that
 $\mathcal{P}$ is divisible by  $\mathcal{L}^2$ in $\mathbb{Z}[X_1, \dots, X_s]$. 
\item[(ii)]
There is a non-constant polynomial $\mathcal{L}$ in $\mathbb{Q}[X_1, \dots, X_s]$ such that
 $\mathcal{P}$ is divisible by  $\mathcal{L}^2$ in $\mathbb{Q}[X_1, \dots, X_s]$. 
\end{itemize}

\begin{remark}
Since $\mathbb{Z}[X_1, \dots, X_s]$ is a subring of $ \mathbb{Q}[X_1, \dots, X_s]$ , (i) implies (ii). The converse  follows readily from Lemma  2.2 of \cite{VZ}.
\end{remark}

\section{Main results}
 
 The following assertion is (contained in)
 Theorem 1.1 of \cite{VZ}.
 
 \begin{theorem}
 If  $\mathcal{P}$  does not satisfy property (a) then the  Euler product $\mathfrak{S}_{\mathcal{P}} \ne 0$ {\sl if and only if} $\mathcal{P}$ is squarefree.
 \end{theorem}
 
 
 {\bf Remark}. After our paper \cite{VZ} was published, Bjorn Poonen  has informed us that 
 the ``{\sl if}'' part of  
Theorem 1.1 had already appeared in  his paper  
   \cite[last paragraph of Section 6, p. 365-366]{P}.
   
\vskip .2cm

Let $P_1, \dots, P_s$ be positive integers and let us put 
$$\mathbf{P}:=(P_1, \dots, P_s) \in \mathbb{Z}^s.$$
 Let us consider the set $S_{\mathcal{P}}(\mathbf{P})$ of $(a_1, \dots, a_s) \in \mathbb{Z}^s$
such that the integer $\mathcal{P}(a_1, \dots, a_s)$ is {\sl squarefree} and
$$|a_j| \le P_j \quad \forall j=1, \dots, s.$$
Let us put
$$N_{\mathcal{P}}(\mathbf{P})=\card(S_{\mathcal{P}}(\mathbf{P})).$$
Clearly, $N_{\mathcal{P}}(\mathbf{P})=0$ if $\mathcal{P}$ enjoys property (a).

The following assertion is (second assertion of) Corollary 1.2 of \cite{VZ}.

\begin{corollary}
If a polynomial $\mathcal{P}$ is not squarefree  then
\begin{equation}
\label{eq:one11}
N_{\mathcal P}(\mathbf P) \ll \frac{P_1\ldots P_s}{\min(P_1,\ldots,P_s)}.
\end{equation}
\end{corollary}

Let us define  
the {\sl upper density} $\mathfrak D_{\mathcal P}$ as
\[
\mathfrak D_{\mathcal P}= \limsup_{\min\{P_1,\ldots,P_s\}\rightarrow\infty} \frac{N_{\mathcal P}(\mathbf P)}{2^sP_1\ldots P_s}
.\]

 The following assertion is Corollary 1.3 of \cite{VZ}.

\begin{corollary}
We have $\mathfrak D_{\mathcal P}\le \mathfrak S_{\mathcal P}$.
 In particular if $\mathfrak D_{\mathcal P}>0$, then $\mathfrak S_{\mathcal P}>0$ and $\mathcal P$ is a squarefree polynomial that does not enjoy property (a).
 \end{corollary}


\begin{thebibliography}{99}



\bibitem{P} B. Poonen, Squarefree values of multivariable polynomials, Duke Math. J. {\bf 118} (2003), 353–373; 	arXiv:math/0203292 [math.NT].



\bibitem{VZ} R. C. Vaughan,  Yu, G. Zarhin, A note on the squarefree density of  polynomials. Mathematika {\bf 70} (2024), no.4, Paper No. e12275, 18 pp.

\end{thebibliography}
\end{document}